\documentclass[12pt]{amsart}
\usepackage{latexsym, amssymb, amsmath}
 \usepackage{eucal}

    \newtheorem{rema}{Remark}[section]
    \newtheorem{propo}[rema]{Proposition}
\newtheorem{algo}[rema]{Algorithm}
   \newtheorem{theo}[rema]{Theorem}
 
   \newtheorem{defi}[rema]{Definition}
    \newtheorem{lemma}[rema]{Lemma}
    \newtheorem{corol}[rema]{Corollary}
     
  \newtheorem{rmk}[rema]{Remark}

	\newcommand{\nno}{\nonumber}

 \newcommand{\pf}{{\it Proof:}\hspace{2ex}}
 
 \newcommand{\epfv}{\hspace{1em}$\Box$\vspace{1em}}

\newcommand{\bC}{{\mathbb C}}

\newcommand{\bN}{{\mathbb N}}
\newcommand{\bT}{{\mathbb T}}

\newcommand{\BQ}{\begin{eqnarray}}
\newcommand{\EQ}{\end{eqnarray}}
\newcommand{\BQn}{\begin{eqnarray*}}
\newcommand{\EQn}{\end{eqnarray*}}
\newcommand{\wtilde}{\widetilde}

\title[A Family of Invariants of Rooted Forests]
{A Family of Invariants of Rooted Forests}

    \author{Wenhua Zhao}      
    \begin{document}

\begin{abstract}
Let $A$ be a  commutative $k$-algebra over a field of $k$ and
 $\Xi$ a linear operator defined on $A$.
We define a family of $A$-valued 
invariants $\Psi$ for finite rooted forests 
by a recurrent algorithm using the operator
$\Xi$ and show that  
the invariant $\Psi$ distinguishes rooted forests 
if (and only if) it distinguishes
rooted  trees $T$, and  
if (and only if) it is {\it finer} than
the quantity $\alpha (T)=|\text{Aut}(T)|$ of 
rooted  trees $T$.
We also consider 
the generating function $U(q)=\sum_{n=1}^\infty U_n q^n$ with 
$U_n =\sum_{T\in \bT_n} \frac 1{\alpha (T)} \Psi (T)$, where 
 $\bT_n$ is
the set of rooted trees with $n$ vertices.
We show that the generating function $U(q)$ satisfies 
the equation $\Xi \exp U(q)= q^{-1} U(q)$. 
Consequently, we get a
recurrent formula for $U_n$ $(n\geq 1)$, namely,
$U_1=\Xi(1)$ and 
$U_n =\Xi S_{n-1}(U_1, U_2, ... , U_{n-1})$ for any
 $n\geq 2$, where $S_n(x_1, x_2, \cdots)$ $(n\in \bN)$ 
are the elementary Schur polynomials.
We also show that the (strict) order polynomials and 
two well known quasi-symmetric function invariants of 
rooted forests are in the family of invariants $\Psi$ 
and derive some 
consequences about these well-known invariants from
our general results on $\Psi$. Finally, we generalize 
the invariant $\Psi$ to labeled planar forests and discuss its
certain relations with 
the Hopf algebra $\mathcal H_{P, R}^D$ in \cite{F} spanned by 
labeled planar forests.
\end{abstract}

   

 \bibliographystyle{alpha}
    \maketitle

\tableofcontents

\renewcommand{\theequation}{\thesection.\arabic{equation}}
\renewcommand{\therema}{\thesection.\arabic{rema}}
\setcounter{equation}{0}
\setcounter{rema}{0}
\setcounter{section}{0}

\section{\bf Introduction}\label{S1}

By a {\it rooted tree} we mean a finite
1-connected graph with one vertex designated as 
its {\it root}. Rooted trees not only form a family of 
important objects in combinatorics,  they are also 
closely related with many other mathematical areas. 
For the connection with the inversion problem and the 
Jacobian problem, see \cite{BCW}, \cite{W}. 
For the connection with D-log and formal flow of analytic maps, 
see \cite{WZ}. For the connection 
with renormalization of quantum field theory, 
see \cite{Kr}, \cite{CK}. 

In this paper, motivated 
by certain properties
 of the (strict) order polynomials encountered in \cite{WZ}, 
we define a family of 
$A$-valued invariants $\Psi$ for rooted forests 
by a recurrent algorithm 
(See Algorithm \ref{algo}) starting with 
an arbitrary commutative
$k$-algebra  $A$ over a field of $k$ and a fixed
linear operator $\Xi$ defined on $A$.  
We show in Proposition \ref{P4.2} and 
Theorem \ref{main2} that the invariant $\Psi$ 
distinguishes rooted forests if (and only if) it distinguishes rooted
trees $T$, and 
if (and only if) it is {\it finer} than 
the quantity $\alpha (T)=|\text{Aut}(T)|$ of rooted
trees $T$. 
In Section \ref{S5},
we consider 
the generating function $U(q)=\sum_{n=1}^\infty U_n q^n$, where 
$U_n =\sum_{T\in \bT_n} \frac 1{\alpha (T)} \Psi (T)$ with $\bT_n$ 
the set of rooted trees with $n$ vertices and $\alpha (T)=|\text{Aut}(T)|$.
We show that the generating function $U(q)$ satisfies the equation
$\Xi \exp U(q)= q^{-1} U(q)$. 
Consequently, we get the recurrent formula $U_1=\Xi(1)$ and
$U_n =\Xi S_{n-1}(U_1, U_2, ... , U_{n-1})$ for any $n\geq 2$, 
where $S_n(x_1, x_2, \cdots )$ are the elementary Schur polynomials.
In Section \ref{S6} and \ref{S7},
we show that,
with properly chosen $A$ and the linear operator $\Xi$,
the (strict) order polynomials of rooted trees and
two families of quasi-symmetric functions for
rooted forests (See (\ref{DefBarPhi}) and 
(\ref{DefPhi}) for the definitions) 
are in the family of invariants $\Psi$ 
defined by Algorithm \ref{algo}.
We also derive some consequences
on these well-known invariants from our 
general results on the invariant $\Psi$. 
Finally, in Section \ref{S8}, 
We generalize 
our invariants to labeled planar forests and discuss certain
relations of our invariants with 
the Hopf algebra $\mathcal H_{P, R}^D$ in \cite{F} spanned by 
labeled planar forests.

The author would like to thank  
Professor John Shareshian, from whom
the author learned Lemma \ref{LL6.3} and 
quasi-symmetric functions for finite posets, 
and Professor David Wight for  many personal 
communications. The author is also very grateful to the referee 
for suggesting to consider the generalization of our invariants to 
labeled planar forests and its relations with 
the Hopf algebra $\mathcal H_{P, R}^D$ spanned by 
labeled planar forests. The last section of this paper 
is the outcome of the efforts to these directions.

\renewcommand{\theequation}{\thesection.\arabic{equation}}
\renewcommand{\therema}{\thesection.\arabic{rema}}
\setcounter{equation}{0}
\setcounter{rema}{0}

\section{\bf Notation and An Operation for Rooted Forests}\label{S2}

{\bf Notation:}
In a rooted tree
there are natural ancestral relations between vertices.  We say a
vertex $w$ is a child of vertex $v$ if the two are connected by an
edge and $w$ lies further from the root than $v$. We define the {\it degree}
of a vertex $v$ of $T$ to be the number of its children. 
A vertex is called a {\it leaf}\/ if it has no
children.  By a rooted forest we mean a disjoint union of finitely many 
rooted trees. When we speak of isomorphisms between rooted forests, we will
always mean root-preserving isomorphisms, i.e. 
the image of a root of a connected component 
which is always a rooted tree  must be a root.

Once for all, we fix the following notation for the rest of this paper.
\begin{enumerate}
\item We let $\mathbb T$ be the set of isomorphism classes of all
rooted trees and $\mathbb F$ the set of isomorphism classes of all
rooted forests.
For $m\ge1$ an integer, we let $\bT_m$ (resp ${\mathbb F}_m$)
the set of isomorphism classes of all rooted trees (resp. forests)
with $m$ vertices.

\item For any rooted tree $T$, we set the following notation:
\begin{itemize}
\item $\text{rt}_T$ denotes the root vertex of $T$.
\item $V(T)$ (resp. $L(T)$) denotes the set of vertices (resp. leaves) of $T$.
\item $v(T)$ (resp. $l(T)$) denotes the
 number of the elements of $V(T)$ (resp. $L(T)$).
\item For $v\in V(T)$ we define the {\it height} of $v$ to be the
number of edges in the (unique) geodesic connecting $v$ to
$\text{rt}_{T}$.
\item $h(T)$ denotes the height of $T$.
\item $\alpha (T)$ denotes the number of the elements of the automorphism
group $\mbox{Aut}(T)$.
\item For $v_{1},\ldots,v_{r}\in V(T)$, we write $T\backslash
\{v_{1},\ldots,v_{r}\}$ for the graph obtained by deleting each of
these vertices and all edges adjacent to these vertices.
\end{itemize}
\item A {\it rooted subtree} of a rooted tree $T$ is defined as a
connected subgraph of $T$ containing $\text{rt}_{T}$, with
$\text{rt}_{T'}=\text{rt}_{T}$.  
\item  We call  the rooted tree with one vertex the {\it singleton}, denoted by 
$\circ$.  
\end{enumerate}

We define the operation $B_+$ for rooted forests as follows. Let 
$S$ be a rooted forest which is disjoint union of rooted trees
$T_i$ 
$(i=1, 2, ..., d)$. We define $B_+(S)=B_+(T_1, ..., T_d)$ to 
be the rooted tree obtained by connecting all roots of  $T_i$ 
$(i=1, 2, ..., d)$ to a single new vertex, which is set to the root of  
the new rooted tree $B_+(T_1, ..., T_d)$. 
If a forest $S$ is the disjoint
union of
$k_1$ copies of rooted tree $T_1$, $k_2$ copies of $T_2$, $\cdots$, 
$k_d$ copies of $T_d$, we also use the notation 
$B_+(T_1^{k_1}, \cdots  T_r^{k_d})$ for the new rooted tree $B_+(S)$.

\begin{lemma}\label{key}
For any $k_i\in \bN$, $T_i\in \bT$ $(i=1, 2, \cdots, d)$ with 
$T_i\ncong T_j$ $(i\neq j)$, we have
\BQ
\alpha (B_+(T_1^{k_1}, \cdots  T_r^{k_d}))=
(k_1)! \cdots (k_d)!
\alpha (T_1)^{k_1} \cdots \alpha (T_d)^{k_d} 
\EQ
\end{lemma}
\pf Set $T=B_+(T_1^{k_1}, \cdots  T_r^{k_d})$ and $R$ the rooted subtree 
of $T$ consisting of the root $\text{rt}_T$ of $T$ 
and all its children in $T$. Let
$\phi: \text{Aut}(T)\to \text{Aut}(R)$ be the restriction map which clearly
is a homomorphism of groups. Let $H\leq \text{Aut(R)}$ 
be the image of $\phi$.
 Since 
$T_i\ncong T_j$ for any $i\neq j$, It is easy to see that
$|H|=k_1!k_2!\cdots k_d!$. Let $K$ be the kernal of $\phi$. 
Note that an element $\alpha \in \text{Aut}(T)$ is in $K$ if and only if
it fixes all the vertices of $R$. Hence  the order
$|K|$ is equals to $\prod_{i=1}^d \alpha (T_i)^{k_i}$. 
Therefore, we have
\BQn
\alpha (T)=|K||H|=(k_1)! \cdots (k_d)!
\alpha (T_1)^{k_1} \cdots \alpha (T_d)^{k_d} 
\EQn
\epfv

\renewcommand{\theequation}{\thesection.\arabic{equation}}
\renewcommand{\therema}{\thesection.\arabic{rema}}
\setcounter{equation}{0}
\setcounter{rema}{0}

\section{\bf A Family of Invariants $\Psi$ for Rooted Forests}\label{S3}

Let $A$ be a commutative $k$-algebra over a field $k$ 
and $\Xi$ an $k$-linear map from $A$ to $A$. Set $a=\Xi(1)$.
We first define an $A$-valued 
invariant $\Psi$ 
for rooted forests by the following algorithm:

\begin{algo}\label{algo} \quad

\begin{enumerate}
\item For any rooted tree $T\in\bT$, we define $\Psi (T)$ as follows. 
\begin{enumerate}
\item[(i)] For each leaf $v$ of $T$,  set $N_v=\Xi(1)=a$.
\item [(ii)] For any other vertex $v$ of $T$, define $N_v$ inductively
starting
from the highest level by setting
$N_v =\Xi (N_{v_1} N_{v_2} \cdots N_{v_k})$, where
$v_j$ $(j=1, 2, \dots, k)$,
 are the distinct
children of $v$.
\item[(iii)] Set $\Psi (T)=N_{\text{rt}_{T}}$.
\end{enumerate}
\item For any rooted forest $S\in \mathbb F$, we set
\BQ\label{PsiForF}
\Psi(S)=\prod_{i=1}^m \Psi (T_i)
\EQ
where $T_i$ $(i=1, 2, \cdots, m)$ are
connected components of $S$.
\end{enumerate}
\end{algo}

From Algorithm \ref{algo}, the following two lemmas are obvious.
\begin{lemma}
$a)$ $B=\{\Psi (S) | S\in {\mathbb F}\}$ 
is a multiplicative subset of $A$, i.e. it is closed under 
the multiplication of the algebra $A$.

$b)$ $\Xi (B) \subset B$.
\end{lemma}

\begin{lemma}\label{Delta}
  Let  $\Gamma$ be an $A$-valued invariant  for rooted forests.
$\Gamma$ can be re-defined and calculated by Algorithm
\ref{algo}  for some $k$-linear map $\Xi$ if and only if 
\begin{enumerate}
\item It satisfies Eq. 
$(\ref{PsiForF})$ for any rooted forest $S\in \mathbb F$.
\item For any
rooted tree $T\in \bT$,
\BQ
  \Gamma(T)=\Xi \prod_{i=1}^d \Gamma (T_i)
\EQ
where $T_i$ $(i=1, 2, ..., d)$ are the connected components 
of $T\backslash \text{rt}_T$.
\end{enumerate}
\end{lemma}

In Section \ref{S6} and \ref{S7}, we will show 
that the strict order polynomials, 
order polynomials as well as two family
 well known quasi-symmetric functions of rooted forests
(See (\ref{DefBarPhi}) and (\ref{DefPhi}) for the definitions) 
are in this family of invariants $\Psi$. 

\renewcommand{\theequation}{\thesection.\arabic{equation}}
\renewcommand{\therema}{\thesection.\arabic{rema}}
\setcounter{equation}{0}
\setcounter{rema}{0}

\section{\bf When the Invariants $\Psi$ Distinguish Rooted Forests}\label{S4}

\begin{defi} We say an invariant $\Psi$ {\it distinguishes} rooted forests 
$(\text{resp. trees})$  if, 
for any $S_1, S_2 \in \mathbb F$ $(\text{resp. $S_1, S_2 \in \mathbb T$})$, 
$\Psi (S_1) =\Psi(S_2)$ if and only if $S_1\simeq S_2$ .
 We say an invariant $\Psi$ is {\it finer} than 
the quantity $\alpha (T)$ if, 
for any $T_1, T_2 \in \bT$,  $\alpha(T_1)\neq \alpha(T_2)$ implies
$\Psi (T_1) \neq \Psi(T_2)$.
\end{defi}

In Combinatorics, it is very desirable to find an invariant
 which can distinguish 
rooted trees or rooted forests. 
For the invariant $\Psi$ defined by Algorithm \ref{algo}, 
we have the following results.

\begin{propo}\label{P4.2}
An $A$-valued invariant $\Psi$ defined by 
Algorithm \ref{algo}  distinguishes rooted forests if $(\text{and only if}\,\,)$ it 
distinguishes rooted trees. 
\end{propo}
\pf Suppose $\Psi$ distinguishes rooted trees. Let $S_1, S_2\in \mathbb F$ with
$\Psi(S_1)=\Psi(S_2)$. We need show that $S_1\simeq S_2$. First, 
by Lemma \ref{Delta},
we have $\Psi (B_+(S_i))=\Xi \Psi(S_i)$ for $i=1, 2$. Hence, 
 $\Psi (B_+(S_1))=\Psi (B_+(S_2))$. Therefore, by our assumption,
we have $B_+(S_1)\simeq B_+(S_2)$, 
which clearly implies $S_1\simeq S_2$.
\epfv

\begin{theo}\label{main2}
An $A$-valued invariant $\Psi$ defined by 
Algorithm \ref{algo}  distinguishes rooted trees if $(\text{and only if}\,\,)$ it
is {\it finer} than $\alpha (T)$.
\end{theo}
\pf Suppose $\Psi$  is finer than
$\alpha (T)$. Let $T_1, T_2\in \bT$ such that 
 $\Psi (T_1) =\Psi(T_2)$. Hence 
$\alpha(T_1)= \alpha(T_2)$. We need show that $T_1\simeq T_2$. 

Suppose $T_1$ and $T_2$ are not isomorphic to each other. 
Let $T=B_+(T_1, T_1)$ and  $T'=B_+(T_1, T_2)$.
 By Lemma \ref{Delta}, we have
\BQn
 \Psi(T) &=& \Xi \left (\Psi(T_1)^2 \right )\\
 \Psi(T') &=&\Xi \left ( \Psi(T_1)\Psi(T_2)\right )
\EQn
Since $\Psi (T_1) =\Psi(T_2)$,  we have
$\Psi(T)=\Psi(T')$. On the other hand, by Lemma \ref{key},
we have
\BQn
\alpha (T)&=&2\alpha(T_1)^2\\
\alpha (T')&=&\alpha(T_1)\alpha(T_2)=\alpha(T_1)^2
\EQn
Therefore, $\alpha (T)\neq \alpha (T')$. So $\Psi$ is not finer than
$\alpha(T)$,  which is a contradiction.
\epfv

\renewcommand{\theequation}{\thesection.\arabic{equation}}
\renewcommand{\therema}{\thesection.\arabic{rema}}
\setcounter{equation}{0}
\setcounter{rema}{0}

\section{\bf A Generating Function for the Invariant $\Psi$ of Rooted Trees}\label{S5}

In this section, we fix an invariant $\Psi$ defined by 
Algorithm \ref{algo} and consider the generating 
function 
\BQ\label{DefForU(q)}
U(q)=\sum_{T\in \bT} \frac 1{\alpha (T)}\Psi (T) q^{v(T)}
\EQ
For any $n\geq 1$, set $U_n=\sum_{T\in \bT_n} \frac 1{\alpha (T)}\Psi (T)$. 
Hence, we have 
$U(q)=\sum_{n=1}^\infty U_n q^n$. We will derive an 
equation satisfied by the generating function $U(q)$, from which
 $\{U_n |n \in \bN \}$ $(n\geq 2)$ 
can be calculated recursively by
using the elementary Schur polynomials.

\begin{theo}\label{main1}
 The generating function  $U(q)$ satisfies the 
 equation
\BQ\label{maineq1}
\Xi \, e^{U(q)}=q^{-1}  U(q)
\EQ
\end{theo}

\pf Consider
\BQn
\Xi \, e^{U (q)}&=& \Xi(1) + \Xi \, \sum_{k=1}^\infty \frac { U^k(q)}{k!}\\
&=& a + \sum_{k=1}^\infty \frac \Xi{k!} \left( \sum_{T\in \bT}
\frac 1 {\alpha (T)} \left ( \Psi (T)q ^{v(T)}\right) \right )^k
\EQn
While, for  the general term of the right hand side of 
the equation above, we have 
\begin{align}
& \frac \Xi {k!} \left (\sum_{T\in \bT}
\frac 1{\alpha (T)}\Psi (T)q ^{v(T)}\right )^k \nno \\
&= \sum_{r=1}^k 
\sum_{\substack{T_1, \cdots, T_r \in \bT\\
T_i\ne T_j (i\ne j)}}
\sum_{\substack {k_1+\cdots +k_r=k\\
k_1, \cdots, k_r\geq 1}} \frac 1{(k_1)! \cdots (k_r)!}
\frac {\Xi \left (\prod_{i=1}^r \Psi (T_i)^{k_i}\right ) }
 {\prod_{i=1}^r \alpha (T_i)^{k_i}} 
q^{\sum_{i=1}^r k_i v(T_i)} \nno \\
\intertext{By Lemma \ref{key} and \ref{Delta}, we have}
&=  \sum_{r=1}^k \sum_{\substack{T_1, \cdots, T_r \in \bT\\
T_i\ne T_j (i\ne j)}}
\sum_{\substack {k_1+\cdots +k_r=k\\
k_1, \cdots, k_r\geq 1}} 
\frac { \Psi (B_+(T_1^{k_1}, \cdots,  T_r^{k_r})}
 {\alpha (B_+(T_1^{k_1}, \cdots,  T_r^{k_r}))} 
q^{v(B_+(T_1^{k_1}, \cdots,  T_r^{k_r}))-1}\nno \\
&= q^{-1}   \sum_{r=1}^k \sum_{T\in T_{r, k+1}} 
\frac { \Psi (T)}
 {\alpha (T)} 
q^{v(T)}\nno
\end{align}
where $T_{r, k+1}$ is the set of equivalence classes of rooted trees with 
$k+1$ vertices and the
degree of the root being exactly $r$.
Therefore, we have
\BQn
\Xi\, e^{U(q)} &=& \Xi(1) +q^{-1} 
\sum_{k=1}^\infty \sum_{r=1}^k \sum_{T\in T_{r, k+1}} 
\frac { \Psi (T)}
 {\alpha (T)} q^{v(T)} \\
&=& \Xi(1)+ q^{-1}  \sum_{\substack {T\in  \bT\\ T\neq \circ}} 
\frac { \Psi (T)}
 {\alpha (T)} q^{v(T)} \\
&=& q^{-1} U(q) \quad \quad (\text{  since $U_1=\Psi(\circ)=\Xi(1)$.}) 
\EQn
\epfv

Recall that the elementary
Schur polynomials $S_n(x)$ $(n\in \bN)$ in $x=(x_1, x_2, \cdots, x_k\cdots )$
are defined by the generating function:
\BQ\label{def4schur}
e^{\sum_{k=1}^\infty x_k q^k}=\sum_{n=0}^\infty S_n(x)q^n
=1+\sum_{n=1}^\infty S_n(x)q^n
\EQ

Note that, if we sign the {\it weight} of the variable $x_k$ 
to be $k$ for any $k\in \bN^+$ and set 
\BQn
\text{wt} (x_{i_1}^{a_1}x_{i_2}^{a_2}\cdots x_{i_d}^{a_d})=\sum_{k=1}^d a_k i_k
\EQn
for any $i_k, a_k \in \bN^+$ $(k=1, 2, \cdots, d)$. 
Then,  for any $n\in \bN$,
 $S_n(x)$ is a polynomial which is homogeneous with respect to weight
with $\text{wt}\, S_n(x)=n$. 
In particular, $S_n(x)$ depends only on the variables $x_i$ 
$(i=1, 2, \cdots, n)$. For more properties of the  elementary
Schur polynomials $S_n(x)$ and their relationship 
with Schur symmetric functions, see \cite{Ka}, 
\cite{M}.

\begin{propo}\label{main-propo} 
For any $n\geq 1$, we have
\BQ
U_1 &=& \Xi\,(1)=a \\
U_n &=& \Xi \, S_{n-1} (U_1, U_2, \cdots, U_{n-1})\label{Xi-recur}
\EQ
\end{propo}
\pf Set $U=(U_1, U_2, \cdots, U_k, \cdots )$.
From (\ref{maineq1}) and (\ref{def4schur}), we have
\BQn
\sum_{n=1}^\infty \Xi\, S_n(U) \, q^n=q^{-1}\sum_{n=2}^\infty 
U_n(t) q^{n-1}
\EQn
By comparing the coefficient of $q^{n-1}$ $(n\geq 2)$, 
we have
\BQn
U_n(t) =  \Xi\,  S_{n-1} (U_1, U_2, \cdots, U_{n-1})
\EQn
Hence, we get (\ref{Xi-recur}).
\epfv

\begin{rmk} One interesting aspect of the invariant $\Psi$ 
and its generating function $U(q)$  is as follows.
From Proposition $\ref{main-propo}$, we see that $U(q)$ 
is the unique solution of Eq. $(\ref{maineq1})$
in the power series algebra $A[[q]]$. Therefore,
any equation of the form
$(\ref{maineq1})$ can be solved  by looking at 
 the invariant $\Psi$ defined by 
Algorithm \ref{algo} for rooted trees and 
its  generating function $U(q)$ defined by $(\ref{DefForU(q)})$.
\end{rmk}

\renewcommand{\theequation}{\thesection.\arabic{equation}}
\renewcommand{\therema}{\thesection.\arabic{rema}}
\setcounter{equation}{0}
\setcounter{rema}{0}

\section{\bf (Strict) Order Polynomials}\label{S6}

Let $T\in\bT$ be a rooted tree. Note that $T$ with the natural partial
order induced from rooted tree structure forms a finite poset 
(partially ordered set), in which the root of $T$ serves the unique 
minimum element. Similarly, any rooted forest also forms a finite poset.
In the rest of this paper, we will always
view rooted forests as finite 
posets in this way. Recall the {\it strict order polynomial} $ \bar\Omega (P)$ 
for a finite poset $P$ is defined to be the unique polynomial 
$ \bar\Omega (P)$ such that $ \bar\Omega (P)(n)$ equals to 
the number of strict order preserving maps $\phi$
from $P$ to the totally ordered set
$[n]=\{1, 2, \cdots, n\}$ for any $n\geq 1$. Here a map $\phi: P\to [n]$
is said to be {\it strict order preserving} if, for any elements $x, y\in P$ 
with $x>y$ in $P$, then $\phi(x)>\phi(y)$ in $[n]$. Also recall that
the {\it order polynomial} $\Omega (P)$ 
for a finite poset $P$ is defined to be the unique polynomial 
$\Omega (P)$ such that $\Omega (P)(n)$ equals to 
the number of order preserving maps $\phi: P\to [n]$
 for any $n\geq 1$. Here a map $\phi: P\to [n]$
is said to be {\it order preserving} if, for any elements $x, y\in P$ 
with $x>y$ in $P$, then $\phi(x)\geq \phi(y)$ in $[n]$. For 
general studies of these two invariants, see \cite{St1}.

In this section, we show that the strict order polynomials 
$ \bar\Omega(T)$ and order polynomials $\Omega(T)$
are both in the family of the invariants $\Psi$  defined 
by Algorithm \ref{algo}. We also derive some consequences 
from our general results on the invariants $\Psi$.

Consider the polynomial ring $\bC [t]$ in one variable $t$ over $\bC$ and
the difference operator $\Delta$, which is defined by
\BQ
\Delta : \bC [t]&\to& \bC [t]\\
f(t)&\to & f(t+1)-f(t)
\EQ

We define the operator $\Delta^{-1}: \bC[t]\to t \bC[t]$ by setting
$\Delta^{-1} (g)$ to be the unique polynomial 
$f\in t\bC [t]$ such that $\Delta (f)=g$ 
for any $g\in \bC[t]$. Note that $\Delta^{-1}: \bC[t]\to t \bC[t]$ is 
well-defined because that, for any polynomial $f\in \bC[t]$, $\Delta (f)=0$
if and only if $f$ is a constant. 

We also define the operator $\nabla$ 
by
\BQn
\nabla: \bC[t]&\to& \bC[t]\\
f(t) &\to&  f(t)-f(t-1)
\EQn
and $\nabla^{-1}$ by setting $\nabla^{-1}(g)$
to be the unique polynomial $f\in t\bC[t]$ such that 
$\nabla (f)=g$ for any $g\in \bC[t]$.

\begin{propo}\label{P6.1}
Let $A=\bC[t]$, then
the strict order polynomials $ \bar \Omega$ 
$(\text{resp. order polynomials } \Omega)$
of rooted forests
can be re-defined and calculated by Algorithm \ref{algo}
with  $\Xi=\Delta^{-1}$
$(\text{resp. } \Xi=\nabla^{-1})$.   
\end{propo}

The proof of the proposition above immediately follows from 
Lemma \ref{Delta}, the fact that 
$ \bar\Omega$ and $\Omega$ also
satisfy Eq. (\ref{PsiForF}) and 
the following lemma
due to John Shareshian. 

\begin{lemma} \label{LL6.3} 
$[\text{J. Shareshian}]$

For any rooted trees $T_i$ $(i=1, 2, \cdots, r)$, we have 
\BQ
\Delta  \bar \Omega (B_+(T_1, T_2, \cdots,  T_r))
&=& \bar \Omega (T_1)  \bar \Omega (T_2) \cdots  \bar \Omega (T_r)\label{E6.3} \\
\nabla \Omega (B_+(T_1,  T_2, \cdots,  T_r))
&=& \Omega (T_1)\Omega (T_2) \cdots  \Omega (T_r)\label{E6.4}
\EQ
\end{lemma}

For the proof of Eq. (\ref{E6.3}),  see the proof of  
Theorem $4.5$ 
in \cite{WZ}. Eq. (\ref{E6.4}) can be proved similarly. 
Actually, Proposion \ref{P6.1} has been proved in \cite{WZ} 
for the strict order polynomials $\bar \Omega$.

Now we consider the corresponding generating functions 
$ \bar U(t, q)=\sum_{T\in \bT}\frac { \bar \Omega (T)}{\alpha (T)}
q^{v(T)}$ and $U(t, q)=\sum_{T\in \bT}\frac {\Omega (T)}{\alpha (T)}
q^{v(T)}$. By Theorem \ref{main1}, we have

\begin{propo}
The generating functions satisfy the equations
\BQ
e^{ \bar U(t, q)} &=& q^{-1}\Delta  \bar U(t, q) \label{Delta-barU}\\
e^{U(t, q)} &=& q^{-1} \nabla U(t, q)\label{Delta-U}
\EQ
\end{propo}

For any $n\geq 1$, we set $ \bar U_n(t)=\sum_{T\in \bT_n} 
\frac { \bar \Omega (T)}{\alpha (T)}$ and
$U_n(t)=\sum_{T\in \bT_n} \frac { \Omega (T)}{\alpha (T)}$. 
By Proposition \ref{main-propo}, we have

\begin{propo}
$a)$ For any $n \geq 2$, we have
\BQ
 \bar U_1 &=& \Delta^{-1}(1)=t \\
 \bar U_n &=&\Delta^{-1} S_{n-1} ( \bar U_1,  \bar U_2, \cdots,  \bar U_{n-1})
\EQ

$b)$ For any $n \geq 2$, we have
\BQ
U_1 &=& \nabla^{-1}(1)=t \\
U_n &=&\nabla^{-1} S_{n-1} (U_1, U_2, \cdots, U_{n-1})
\EQ
\end{propo}

 Set $u(q)=U(t, 1)$ and
write $u(q)=\sum_{n=1}^\infty u_n q^n$. Since 
 $\Omega (T)(1)=1$ and $\Omega (T)(0)=0$ for any rooted tree $T$, we see that
$u_n=\sum_{T\in \bT_n} \frac 1{\alpha(T)}$ and $U(0, q)=0$.  Therefore,
\BQn
(\nabla U(t, q))(1)=U(1, q)-U(0, q)=u(q)
\EQn
Combining with  Eq. (\ref{Delta-U}), we see that the generating function 
$u(q)$ satisfies the equation
\BQ\label{special-Delta}
e^{u(q)}=q^{-1}u(q)
\EQ
But, on the other hand, it is well known that
there is another generating function related with 
rooted trees satisfying  Eq. (\ref{special-Delta}) which is defined as
follows. Let $r(n)$ be the number of rooted trees on the labeled set 
$[n]=\{1, 2, \cdots, n\}$. Let $R(q)=\sum_{n\geq 1}\frac { r(n)}{n!} q^n$. 
Then, by Proposition $5.3.1$ in \cite{St2},  $R(q)$ also satisfies Eq.
(\ref{special-Delta}) and by Proposition $5.3.2$ in \cite{St2}, we know that
$r(n)=n^{n-1}$.
Therefore, we have 

\begin{corol}
$u(q)=R(q)$. In particular, 
for any $n\geq 1$, we have the identities
\BQ
 n!\sum_{T\in \bT_n} \frac 1{\alpha(T)} &=& r(n)\\
\sum_{T\in \bT_n} \frac 1{\alpha(T)} &=& \frac{n^{n-1}}{n!}
\EQ
\end{corol}

For the corollary above, we see that Eq. (\ref{Delta-U})
can be viewed as a natural generalization of Eq. (\ref{special-Delta}).

\renewcommand{\theequation}{\thesection.\arabic{equation}}
\renewcommand{\therema}{\thesection.\arabic{rema}}
\setcounter{equation}{0}
\setcounter{rema}{0}

\section{\bf Two Quasi-Symmetric Function Invariants for Rooted Forests}\label{S7}

Let us first recall the following well known
quasi-symmetric functions 
$ \bar K (P)$ and $K (P)$ defined 
in \cite{St2} for finite posets $P$. For more general studies on 
quasi-symmetric functions, 
see \cite{G}, \cite{T}, \cite{MR} and \cite{St2}.

Let $x=(x_1, x_2, \cdots, )$ be a sequence of commutative variables and 
$\bC[[x]]$ the formal power series algebra in $x_k$ $(k\geq 1)$ over $\bC$.
For any finite poset 
$P$ and any map $\sigma : P\to \bN^+$ of sets, we set
$x^\sigma:=\prod_{i=1}^\infty x_i^{|\sigma^{-1}(i)|}$ and
define
\BQ\label{DefBarPhi}
 \bar K (P)(x)=\sum_{\sigma} x^\sigma 
\EQ
where the sum runs over the set of all 
strict order preserving maps
$\sigma : P\to \bN^+$. 
Similarly, we define
\BQ\label{DefPhi}
 K (P)(x)=\sum_{\sigma} x^\sigma 
\EQ
where the sum runs over the set of all 
 order preserving maps
$\sigma : P\to \bN^+$. 
 Note that $ \bar K (P)(x)$ and $K(P)(x)$
are always in $\bC [[x]]$ and satisfy Eq. (\ref{PsiForF}) 
for rooted forests.

Recall that an element $f\in \bC[[x]]$ is said to be 
{\it quasi-symmetric} if the degree of $f$ is bounded,
and for any $a_1, a_2, \cdots, a_k \in \bN^+$, $i_1<i_2<\cdots <i_k$ and
 $j_1<j_2<\cdots <j_k$, 
the coefficient of the monomial 
$x_{i_1}^{a_1}x_{i_2}^{a_2}\cdots x_{i_k}^{a_k}$ is 
always same as the coefficient of the monomial 
$x_{j_1}^{a_1}x_{j_2}^{a_2}\cdots x_{j_k}^{a_k}$. From the 
definitions (\ref{DefBarPhi}) and (\ref{DefPhi}), it is easy to check 
that, for any finite
poset $P$,  $ \bar K (P)$ and $K (P)$ 
are quasi-symmetric.

In this section, we will show that the quasi-symmetric functions
$ \bar K$ and $K$ for rooted forests
are also in the family of the invariants $\Psi$ 
defined by Algorithm \ref{algo}.

We define the shift operator $S: \bC[[x]]\to \bC[[x]]$ by 
first setting 
\BQn
S (1)&=&1\\
S (x_m)&=&x_{m+1}
\EQn
and then extending it to $\bC[[x]]$
to be the unique $\bC$-algebra
 homomorphism from $\bC[[x]]$ to $\bC[[x]]$.
For any $m\in \bN^+$, we denote by the abusing notation $x_m$ the
$\bC$-linear map from $\bC[[x]]$ to $\bC[[x]]$ induced by the multiplication
by $x_m$. 

The following lemma follows immediately from the definition of the linear
operator $S$.

\begin{lemma}\label{L7.1}
As the linear maps from  $\bC[[x]]$ to $\bC[[x]]$, 
$x_m S^k=S^k x_{m-k}$ for any $k, m\in \bN^+$ with $k<m$.
\end{lemma}

We define the linear maps $\bar \Lambda$ and $\Lambda$ from
$\bC[[x]]$ to $\bC[[x]]$
 by setting
\BQ
 \bar \Lambda &=& \sum_{k=1}^\infty x_k S^k=
(\sum_{k=1}^\infty S^k) \, x_1 S \label{Bar-Lambda} \\
 \Lambda&=&\sum_{k=1}^\infty x_k S^{k-1}=
(\sum_{k=1}^\infty S^k) \, x_1  \label{Lambda}
\EQ
where the 
last equalities of the 
equations above follow from Lemma \ref{L7.1}.
It is easy to see that $\bar \Lambda$ and $\Lambda$ are well defined. 

\begin{lemma}\label{L7.2}
$a)$ The linear maps $ \bar \Lambda$ and $\Lambda $ from $\bC[[x]]$ to $\bC[[x]]$
are injective.

$b)$ 
\BQ
 \bar\Lambda (1)=\Lambda (1)=\sum_{k=1}^\infty x_k
\EQ
\end{lemma}
\pf $b)$ follows immediately from Eq.
(\ref{Bar-Lambda}) and (\ref{Lambda}).

To prove  $a)$,
let $f \in \bC[[x]]$ such that $ \bar\Lambda f=0$.
By Eq. (\ref{Bar-Lambda}), we have $(1-S) \bar \Lambda =x_1 S$. Hence,
 $x_1 S f=0$ and $Sf=0$. Therefore, 
we must have $f=0$. The injectivity of $\Lambda$ can be proved similarly. 
\epfv

\begin{lemma}
For any rooted tree $T$, we have
\BQ
 \bar K(T)&=& \bar \Lambda \prod_{i=1}^d  \bar K(T_i)\label{Delta-Phi}\\
 K(T)&=& \Lambda \prod_{i=1}^d  K(T_i)\label{Delta-BarPhi}
\EQ
where $T_i$ $(i=1, 2, \cdots, T_d)$ are the connected components
of $T\backslash \text{rt}_T$.
\end{lemma}

\pf Here we only prove Eq. (\ref{Delta-Phi}). For Eq. 
(\ref{Delta-BarPhi}), the ideas of the proof are similar.

Let $W $ be the set of all strict order 
preserving maps $\sigma: P\to \bN^+$ and
 $W_k$ $(k\geq 1)$ the set  
 of $\sigma \in W$ such that 
$\sigma (\text{rt}_T)=k$. Clearly, $W$ equals to the disjoint union
of  $W_k$ $(k\geq 1)$. By the definition of $ \bar K$, 
see (\ref{DefBarPhi}),  we  see that
$\sum_{\sigma\in W_k} x^\sigma \in  \bC[[x_k, x_{k+1}, \cdots, ]]$. Since 
$\bar K$  satisfies Eq. (\ref{PsiForF}) 
for rooted forests, we have
\BQ
\sum_{\sigma\in W_k} x^\sigma = x_k S^k \bar K (T\backslash \text{rt}_T)=
x_k S^k \prod_{i=1}^d   \bar K(T_i)
\EQ
Therefore, 
\begin{align*}
 \bar K(T)&=\sum_{k=1}^\infty \sum_{\sigma\in W_k} x^\sigma\\
&=\sum_{k=1}^\infty x_k S^k \prod_{i=1}^d K(T_i)\\
&= \Lambda (\prod_{i=1}^d K(T_i))
\end{align*}
\epfv
 
From the lemma above and Lemma \ref{Delta} and the fact that
$\bar K$ and $K$  satisfy Eq. (\ref{PsiForF}) 
for rooted forests, we immediately have

\begin{propo}\label{P7.4}
The quasi-symmetric functions $ \bar K$ $(\text{resp. $K$})$
for rooted forests
can be re-defined and calculated by Algorithm \ref{algo}
with $A=\bC[[x]]$ and
$\Xi= \bar\Lambda$ $(\text{resp. $\Xi=\Lambda$})$.
\end{propo}

Now we consider the generating functions
\BQn
 \bar Q(x, q)&=&\sum_{T\in \bT} \frac{ \bar K(P)(x)}{\alpha (T)}q^{v(T)}
=\sum_{n=1}^\infty Q_n(x)q^n\\
Q(x, q)&=&\sum_{T\in \bT_n} \frac{K(P)(x)}{\alpha (T)}q^{v(T)}
=\sum_{n=1}^\infty Q_n(x)q^n
\EQn
where $ \bar Q_n(x)=\sum_{T\in \bT_n} \frac{ \bar K(P)(x)}{\alpha (T)}$ and
$Q_n(x)=\sum_{T\in \bT_n} \frac{K(P)(x)}{\alpha (T)}$ for any $n\geq 1$.
By Theorem \ref{main1}, Lemma \ref{L7.2} and Proposition \ref{P7.4},  we have

\begin{propo} $a)$ The generating functions
$Q (x, q)$ and  $Q (x, q)$ satisfy the equations 
\BQ
 \bar \Lambda e^{ \bar Q(x, t)}&=&q^{-1} Q(x, t)\\
 \Lambda e^{Q(x, t)}&=&q^{-1} Q(x, t)
\EQ

$b)$ Consequently,  we have the
recurrent formula for $ \bar Q_n(x)$ and $Q_n(x)$ $(n\in \bN^+)$
\BQ
 \bar Q_1(x)&=&\sum_{k=1}^\infty x_k\\
 \bar Q_n(x)&=& \bar \Lambda \left (S_{n-1}( \bar Q_1(x), 
 \bar Q_2(x), \cdots, \bar Q_{n-1}(x))\right )
\EQ
and
\BQ
Q_1(x)&=&\sum_{k=1}^\infty x_k\\
Q_n(x, t)&=&\Lambda \left (S_{n-1}(Q_1(x), Q_2(x), \cdots, Q_{n-1}(x))\right )
\EQ
respectively.
\end{propo}

One natural question one may ask is whether or not
the invariants $ \bar\Omega(T)$, 
$\Omega(T)$, $ \bar K(T)$ and $K(T)$ distinguish rooted forests.
 The answers for the strict 
order polynomials $ \bar\Omega $ and order polynomial 
$\Omega$ are well known to be negative. 
(See, for example,  Exercise $3.60$ in \cite{St1}).
For the quasi-symmetric polynomial invariants
$ \bar K$ and $K$, the answers seem to be positive, 
but we do not know
any proof in literature.

One remark is that the invariant $\Psi$ defined by Algorithm 
\ref{algo} can also be extended to the set of finite posets 
by a more general recurrent procedure. 
This will be done in the appearing
paper \cite{SWZ}. But for the corresponding 
generating function
\BQ
V(q)=\sum_P \frac {\Psi(P)}{\alpha(P)}q^{v(P)}
\EQ
where the sum runs over the set of all finite posets $P$, 
it is not clear what the generalization of Eq. (\ref{maineq1})
satisfied by $V(q)$ should be. 
This is unknown even for the case of (strict) order polynomials.

\renewcommand{\theequation}{\thesection.\arabic{equation}}
\renewcommand{\therema}{\thesection.\arabic{rema}}
\setcounter{equation}{0}
\setcounter{rema}{0}

\section{\bf Generalization to Labeled Planar Forests} \label{S8}

In this section, we first generalize the construction of the invariant
$\Psi$ defined by Algorithm \ref{algo} for rooted forests 
to {\it labeled planar forests} and then consider 
its certain relationships with the Hopf algebra 
$\mathcal H_{P, R}^D$ in \cite{F} 
spanned by labeled planar forests. 

Once for all, we fix a non-empty 
finite or countable set $D$. By a labeled planar rooted tree $T$, we always 
mean in this section a rooted tree $T$ such that
each vartex of $T$ is assigned a unique element of $D$ and set of all 
children of any single vertex of $T$ is an ordered set. A labeled planar rooted
forest $F$ is an ordered set of finitely many labeled planar rooted trees.
We let $\mathbb T_{P, R}^D$ denote the set of 
all labeled planar rooted trees
and $\mathbb F_{P, R}^D$ 
the set of all labeled planar rooted forests.
For any labeled planar rooted forest $F=T_1T_2\cdots T_d$,
with $T_i\in \mathbb T_{P, R}^D$ 
$(1\leq i\leq d)$ and $\alpha \in D$,
we define 
$B_+^\alpha (F)=B_+^\alpha (T_1T_2\cdots T_d)$ to be the labeled planar 
rooted tree obtained by connecting the root of each 
$T_i$ to a $\alpha$-labeled vertex $v$ by an edge and set 
the new vertex $v$ to be the root of this new labeled planar 
rooted tree.

We also fix an associative (not necessarily commutative) algebra  $A$ 
over a field $k$ and $\{ \Xi_\alpha |\alpha \in D\}$ a sequence of 
linear operators of $A$. Now we define an $A$-valued invariant 
$\Psi (F)$ for labeled planar forests $F$ 
by the following algorithm.

\begin{algo}\label{algo-2}

\begin{enumerate}
\item For any labeled planar rooted tree $T\in\bT$, 
we define $\Psi (T)$ as follows. 
\begin{enumerate}
\item[(i)] For each $\alpha$-labeled 
leaf $v$ of $T$,  set $N_v=\Xi_\alpha (1)$.

\item [(ii)] For any other vertex $v$ of $T$, define $N_v$ inductively
starting
from the highest level by setting
$N_v =\Xi_\alpha (N_{v_1} N_{v_2} \cdots N_{v_k})$, where
$\alpha$ is the label of $v$ and
$(v_1, v_2, \dots, v_k)$
are the ordered children of $v$.

\item[(iii)] Set $\Psi (T)=N_{\text{rt}_{T}}$.
\end{enumerate}
\item For any labeled planar rooted forest $F=T_1 T_2\cdots T_m$, 
where $T_i$ $(i=1, 2, \cdots, m)$ are
connected components of $F$, 
we set
\BQ\label{M-PsiForF}
\Psi(F)=\Psi (T_1)\Psi (T_2) \cdots \Psi (T_m)
\EQ
\end{enumerate}
\end{algo}

Note that the order in the product in Eq. (\ref{M-PsiForF}) must be same
as the one in the expression $F=T_1T_2\cdots T_m$.  

From Algorithm \ref{algo-2}, the following lemma is obvious.

\begin{lemma}\label{M-Delta}
Let  $\Gamma$ be an $A$-valued invariant 
for labeled planar rooted forests ${\mathbb F}_P^D$. Then
$\Gamma$ can be re-defined and calculated by Algorithm
\ref{algo} for some $k$-linear map $\Xi$ if and only if 
\begin{enumerate}

\item It satisfies Eq. $(\ref{M-PsiForF})$ for any 
$F\in \mathbb F_{P, R}^D$.

\item For any $T\in \bT_{P, R}^D$ with $T=B_+^\alpha (T_1T_2\cdots T_d)$, 
we have
\BQ\label{E8.1.2}
 \Gamma(T)=\Xi_\alpha  (\Gamma (T_1)\Gamma (T_2)\cdots \Gamma (T_d))
\EQ
\end{enumerate}
\end{lemma}

\begin{rmk}\label{R8.1.3}
Let $\mathcal H_{P, R}^D$ be the vector spaces spanned 
by labeled planar forests. In \cite{F}, a Hopf algebra 
structure in $\mathcal H_{P, R}^D$ is given, which is a 
labeled planar version of Kreimer's 
Hopf algebra $(\text{See \cite{Kr} and \cite{CK}})$ spanned by 
rooted forests. The product of the Hopf algebra 
$\mathcal H_{P, R}^D$
is given by the ordered disjoint union operation.
We extend the map $\Gamma$ defined by Algorithm $\ref{algo-2}$ 
to $\mathcal H_{P, R}^D$ linearly and still denote it by $\Gamma$. Then 
it is easy to see that condition $(1)$ in the lemma above 
is equivalent to saying that the map $\Gamma$ is a homomorphism 
of algebras from to  $\mathcal H_{P, R}^D$ to $A$, while condition $(2)$
is equivalent to the following equation.
\BQ\label{Cocycle-Eq}
\Gamma\circ B_+^\alpha =\Xi_\alpha \circ \Gamma
\EQ
\end{rmk}

Now let us consider the corresponding generating functions 
$U(q)$ for the 
invariants of $\Psi$ defined by Algorithm \ref{algo-2}. 

First, for each $\alpha \in D$, we set
\BQ
U_\alpha (q)=
\sum_{T\in \bT_{P, R, \alpha}^D} \Psi (T) q^{v(T)}
\EQ
where $\bT_{P, R, \alpha}^D$ is the set of all labeled planar 
rooted trees with $\alpha$-labeled roots. 
We also set
\begin{align}
U(q)&=\sum_{\alpha\in D}U_\alpha (q) \label{E8.1.5}\\
\Xi &=\sum_{\alpha\in D} \Xi_\alpha
\end{align}

\begin{theo}\label{T8.1.4}
The generating functions $U_\alpha (q)$ $(\alpha \in D)$ and $U(q)$
satisfy the following equations.
\begin{align}
\Xi_\alpha \, \frac 1{1-U(q)}& =q^{-1}  U_\alpha (q)\label{E8.1.7} \\
\Xi \, \frac 1{1-U(q)} & =q^{-1}  U (q)\label{EqForPlanar}
\end{align}
\end{theo}

First, note that the second equation follows from the first 
one by taking sum over the set $D$. 
The proof of the first equation is parallel to the proof 
of Eq. (\ref{maineq1})
in Theorem \ref{main1} but a little easier, since 
automorphism groups of planar labeled rooted trees are trivial. 
So we omit the proof here.

\begin{rmk}\label{R8.1.5}

$(1)$  Note that, when $|D|=1$, $\mathbb F_{P, R}^D$ is same as the set 
$\mathbb F_{P, R}$ of unlabeled planar rooted forests. Hence 
Algorithm $\ref{algo-2}$ gives an invariant 
for planar rooted forests in this case.
Since the solution of Eq. $(\ref{EqForPlanar})$ in $A[[q]]$ is unique, 
any equation of the form Eq.
$(\ref{EqForPlanar})$ can be solved by looking at 
the invariant $\Psi$ defined by 
Algorithm $\ref{algo-2}$ for planar rooted trees and 
its generating function $U(q)$ defined by 
Eq. $(\ref{E8.1.5})$.

$(2)$ When $|D|=1$ and the algebra $A$ is commutative, for any 
planar rooted forest $F$, the invariant $\Psi (F)$ defined 
by Algorithm $\ref{algo-2}$ coincides with the one defined by 
Algorithm $\ref{algo}$ for the underlying rooted forest of $F$, 
which is obtained by simply ignoring the planar structure of $F$.
\end{rmk}

Next we discuss certain relationships of the invariants $\Psi$ defined by
Algorithm \ref{algo-2} with the Hopf algebra $\mathcal H_{P, R}^D$
defined and studied in \cite{F}. Even though, the links present here have no
obvious logical implication one way or the other, they provide a 
new point of view to the invariants $\Psi$ defined by Algorithm \ref{algo-2}.
Besides the Hopf algebra $\mathcal H_{P, R}^D$ 
and its certain universal property studied in 
\cite{F}, we also need the Hopf algebra structure  
defined in \cite{F} on the tensor algebra $T(V)$ 
for any vector space $V$.
Since definitions of various operations of the Hopf algebras 
$\mathcal H_{P, R}^D$ and $T(V)$ are quite
involved,  we will follow the notation in \cite{F} closely
and quote necessary results directly from \cite{F}.
We refer readers to \cite{F} and references there 
for more details.

First, let us assume that our fixed associate algebra 
$(A, m, \eta)$ also has a co-algebra structure 
with which it forms
a bi-algebra $(A, m, \eta, \Delta, \epsilon)$. We further assume that 
the linear operators $\Xi_\alpha$ $(\alpha \in D)$ 
are 1-cocycles, i.e. they satisfy the following equation
\begin{align}
\Delta \circ \Xi_\alpha=\Xi_\alpha \otimes 1
+(\text{id} \otimes \Xi_\alpha)\circ \Delta .
\end{align}

By the universal property of the Hopf algebra 
$\mathcal H_{P, R}^D$ given in Theorem $24$ in \cite{F},
there exists a unique homomorphism of bi-algebras 
$\varphi: \mathcal H_{P, R}^D \to A$ such that
\BQ\label{EE8.1.11}
\varphi \circ B_\alpha^+ = L_\alpha \circ \varphi . 
\EQ

Note that the map $\varphi: \mathcal H_{P, R}^D \to A$
gives an $A$-valued invariant for labeled planar forests.

\begin{propo}\label{PP8.1.6}
The $A$-valued invariant
$\varphi (F)$ defined above belongs to the family of invariants of 
labeled planar forests defined by Algorithm $\ref{algo-2}$ with 
the linear operators $L_\alpha$ $(\alpha \in D)$.
\end{propo}

In other words, in this special situation, 
 the invariant 
$\Psi$ by Algorithm \ref{algo-2} coincides with the 
unique map $\varphi$ guaranteed by the universal 
property of the Hopf algebra
$\mathcal H_{P, R}^D$. \newline

\pf 
Since the homomorphism $\varphi$ preserves the algebra products 
and satisfies Eq. (\ref{Cocycle-Eq}),  the proposition follows 
immediately from Remark $\ref{R8.1.3}$ and Lemma $\ref{M-Delta}$.
\epfv

One remark is that Algorithm \ref{algo-2} does not depends on whether 
the algebra $A$ has a bi-algebra structure.
It only depends on the associate
algebra structure of $A$. But, on the other hand, it is 
shown in \cite{F} that the tensor algebra $T(V)$ 
of any vector space $V$ 
has a Hopf algebra structure. In particular, we have a Hopf 
algebra sturcture on the tensor algebra $T(A)$. 
Next we  
show that, by using the linear operators $\Xi_\alpha$ $(\alpha \in D)$
and the associate algebra stucture of $A$, we can construct 
a family of 1-cocycles $L_\alpha$ $(\alpha \in D)$ of 
the Hopf algebra $T(A)$. Therefore,
the corresponding unique map $\varphi: \mathcal H_{P, R}^D \to T(A)$ 
does give us a family of $T(A)$-valued invariants for labeled planar 
forests. 

First, for any $\alpha \in D$,  we define a linear map
from $\wtilde \Xi_\alpha: T(A)\to A$ by setting
\begin{align*}
\wtilde \Xi_\alpha: k &\to A\\
a &\to a \, \Xi (1_A) 
\end{align*}
and, for any $n\geq 1$,  
\begin{align*}
\wtilde \Xi_\alpha: A^{\otimes n} &\to A\\
v_1\otimes v_2\cdots \otimes v_n &\to  \Xi_\alpha ( 
v_1 \cdot v_2\cdots  v_n)
\end{align*}
and extend it linearly to $T(A)$. Note that
here we use $1_A$ for the identity element of 
the algebra $A$ to distinguish the identity element $1_k$ 
in the ground field $k$.

Next we define a sequence linear maps 
$\{L_\alpha : T(A)\to T(A) |\alpha \in D\}$ 
by setting 
\begin{align}
L_\alpha (a)&=a \, \Xi(1_A) \quad \text{for any $a\in k$ and } \label{E8.1.9}\\
L_\alpha & (v_1\otimes v_2\cdots \otimes v_n)\label{E8.1.10}\\
&= \sum_{j=1}^{n-1}
v_1\otimes v_2\cdots \otimes v_j \otimes
\wtilde \Xi_\alpha ( v_{j+1} \otimes \cdots\otimes v_n) \nno \\
& \quad +\wtilde \Xi_\alpha ( v_{1} \otimes \cdots\otimes v_n)+
v_1\otimes v_2\cdots \otimes v_n\otimes \wtilde \Xi_\alpha (1)\nno \\
&= \sum_{j=1}^{n-1}
v_1\otimes v_2\cdots \otimes v_j \otimes
 \Xi_\alpha ( v_{j+1}\cdot v_{j+2} \cdots v_n) \nno \\
& \quad + \Xi_\alpha ( v_{1}\cdot v_{2} \cdots v_n)+
v_1\otimes v_2\cdots \otimes v_n\otimes  \Xi_\alpha (1)\nno
\end{align} 
and extend it linearly to $T(A)$.

By Proposition $72$ in \cite{F}, 
the linear maps $L_\alpha : T(A)\to T(A)$ are 1-cocycles 
of the Hopf algebra $T(A)$.
By the universal property of the Hopf algebra 
$\mathcal H_{P, R}^D$ given in Theorem $24$ in \cite{F},
there exists a unique homomorphism of Hopf algebras 
$\varphi: \mathcal H_{P, R}^D \to T(A)$ such that
\BQ\label{E8.1.11}
\varphi \circ B_\alpha^+ =L_\alpha \circ \varphi . 
\EQ

Note that the map $\varphi: \mathcal H_{P, R}^D \to T(A)$ gives  
a $T(A)$-valued invariant for labeled planar forests, which, 
by Proposition \ref{PP8.1.6}, is same as the one
defined by Algorithm $\ref{algo-2}$ with 
the linear operators $L_\alpha$ $(\alpha \in D)$.

{\small \sc Department of Mathematics, Washington University in St.
Louis,
St. Louis, MO 63130-4899. }

{\em E-mail}: 
zhao@math.wustl.edu.

\end{document}